\documentclass[a4paper,11pt]{amsart}
\usepackage{graphicx}
\usepackage{mathptmx}
\usepackage{amsmath}
\usepackage{amssymb}
\usepackage{enumitem}
\usepackage{xcolor}

\newmuskip\pFqmuskip

\newcommand*\pFq[6][8]{%
  \begingroup 
  \pFqmuskip=#1mu\relax
  \mathcode`=\string"8000
  \begingroup\lccode`\~=`\,
  \lowercase{\endgroup\let~}\pFqcomma
  F^{#2}_{#3}{\left(\genfrac..{0pt}{}{#4}{#5}\bigg|#6\right)}%
  \endgroup
}
\newcommand{\pFqcomma}{\mskip\pFqmuskip}

\newtheorem{theorem}{Theorem}

\newtheorem{remark}[theorem]{Remark}

\begin{document}

\title[]{generalized degenerate Stirling numbers arising from degenerate Boson normal ordering }

\author{Taekyun Kim $^{1}$}
\address{Department of Mathematics, Kwangwoon University, Seoul 139-701, Republic of Korea}
\email{tkkim@kw.ac.kr}

\author{Dae San Kim $^{2}$}
\address{Department of Mathematics, Sogang University, Seoul 121-742, Republic of Korea}
\email{dskim@sogang.ac.kr}

\author{Hye Kyung Kim $^{3,*}$}
\address{Department of Mathematics Education, Daegu Catholic University, Gyeongsan 38430, Republic of Korea}
\email{hkkim@cu.ac.kr}

\thanks{* is corresponding author}

\begin{abstract}
It is remarkable that, in recent years, intensive studies have been done for degenerate versions of many special polynomials and numbers and have yielded many interesting results. The aim of this paper is to study the generalized degenerate $(r,s)$-Stirling numbers of the second and their natural extensions to polynomials, namely the generalized degenerate $(r,s)$-Bell polynomials, arising from certain `degenerate boson normal ordering.' We derive some properties, explicit expressions and generating functions for those numbers and polynomials. The generalized degenerate $(r,s)$-Stirling numbers of the second and the degenerate boson normal ordering are respectively degenerate versions of the generalized $(r,s)$-Stirling numbers of the second and the boson normal ordering studied earlier by Blasiak-Person-Solomon.
\end{abstract}

 \maketitle

\markboth{\centerline{\scriptsize Generalized degenerate Stirling numbers arising from degenerate Boson normal ordering}}
{\centerline{\scriptsize T. Kim, D. S. Kim, and H. K. Kim}}


\section{Introduction}
The Stirling number of the second kind $S_2(n,k)$ is the number of ways to partition a set of $n$ objects into $k$ nonempty subsets (see \eqref{eq0}). The Stirling numbers of the second kind arise in various different contexts and have numerous applications, for example to enumerative combinatorics and quantum mechanics. In recent years, intensive explorations have been done for degenerate versions of many special numbers and polynomials, which was initiated by Carlitz in his work on degenerate Bernoulli and degenerate Euler polynomials. They have been studied by using such tools as combinatorial methods, generating functions, $p$-adic analysis, umbral calculus techniques, probability theory, mathematical physics, operator theory, special functions, analytic number theory and 
differential equations. The degenerate Stirling numbers of the second kind $S_{2,\lambda}(n,k)$ appear naturally when we replace the power $x^n$ by the generalized falling factorial polynomial $(x)_{n,\lambda}$ in the defining equation of $S_2(n,k)$ (see \eqref{eq0}, \eqref{eq02}). It turns out that they appear very frequently when we study degenerate versions of many special polynomials and numbers. 

The normal ordering of $(a^{\dag}a)^{n}$ in terms of boson operators $a$ and $a^{\dag}$ can be written in the form
\begin{equation}\label{eq -2}
(a^{\dag}a)^{n}=\sum_{l=0}^{n}S_{2}(n,l)(a^{\dag})^{l}a^{l}.
\end{equation}
The normal ordering of $(a^{\dag}a)_{n,\lambda}$ in terms of boson operators $a$ and $a^{\dag}$, which is a degenerate version of \eqref{eq -2}, can be written in the form
\begin{equation*}
(a^{\dag}a)_{n,\lambda}=\sum_{l=0}^{n}S_{2,\lambda}(k,l)(a^{\dag})^{l}a^{l},\label{1}
\end{equation*}
where the generalized falling factorials $(x)_{n,\lambda}$ are given by \eqref{eq01-1}.\par
Let $r,s$ be positive integers with $r\geq s$, and let $n$ be any positive integer. Blasiak-Person-Solomon introduced the generalized $(r,s)$-Stirling numbers of the second kind $S_{r,s}(n,k)$, which boil down to $S_{1,1}(n,k)=S_{2}(n,k)$, for $r=s=1$, by considering the boson normal ordering of $((a^{\dag})^r a^s)^n$:
\begin{equation}\label{eq -1}
\begin{split}
((a^{\dag})^r a^s)^n=(a^{\dag})^{n(r-s)}\sum_{k=s}^{ns}S_{r,s}(n,k)(a^{\dag})^k a^k.
\end{split}
\end{equation}
In this paper, we introduce the generalized degenerate $(r,s)$-Stirling numbers of the second kind, which are degenerate versions of the generalized $(r,s)$-Stirling numbers of the second kind, by considering a degenerate version of \eqref{eq -1}, namely the boson normal ordering of $ \prod_{k=0}^{n-1}\big((a^{\dag})^{r}a^s -k\lambda (a^{\dag})^{r-s}\big)$:
\begin{equation*}
\begin{split}
\prod_{k=0}^{n-1}\Big[(a^{\dag})^{r-s}\Big((a^{\dag})^s a^s -k\lambda\Big)\Big]
&=\prod_{k=0}^{n-1}\Big((a^{\dag})^{r}a^s -k\lambda (a^{\dag})^{r-s}\Big)\\
&=(a^{\dag})^{n(r-s)}\sum_{k=0}^{ns}S_\lambda^{(r,s)}(n,k)(a^{\dag})^k a^k.
\end{split}
\end{equation*} \par
The aim of this paper is to derive some properties, explicit expressions and generating functions for the generalized degenerate $(r,s)$-Stirling numbers of the second kind and their natural extensions to polynomials, namely the generalized degenerate $(r,s)$-Bell polynomials. The novelty of this paper is that the generalized degenerate $(r,s)$-Stirling numbers of the second kind are introduced in a natural manner by considering the `degenerate boson normal ordering.' We think that these new numbers will play an important role in the study of various degenerate versions of many special polynomials and numbers. \par
In more detail, the outline of this paper is as follows. We derive several expressions for the generalized degenerate $(r,s)$-Bell polynomials $\phi_{n,\lambda}^{(r,s)}(x)$ (see \eqref{eq15}) in Theorem 2 and the generalized degenerate $(r,s)$-Bell numbers $\phi_{n,\lambda}^{(r,s)}=\phi_{n,\lambda}^{(r,s)}(1)$ in Theorems 2 and 3. We obtain several expressions for the generalized degenerate $(r,s)$-Stirling numbers of the second kind in Theorems 4-6. $\phi_{n,\lambda}^{(r,r)}(|z|^2)$ and its generating function $\sum_{n=0}^{\infty}\phi_{n,\lambda}^{(r,r)}(|z|^2)\frac{t^n}{n!}$ are expressed in terms of bra-ket notation respectively in Theorem 7 and Theorem 8. We deduce the generating function $\sum_{n=0}^{\infty}\phi_{n,\lambda}^{(r)}(|z|^2)\frac{t^n}{n!}$ of the degenerate $r$-Bell polynomials $\phi_{n,\lambda}^{(r)}(x)$, which are different from $\phi_{n,\lambda}^{(r,r)}(x)$ and natural extension to polynomials of the degenerate $r$-Stirling numbers of the second kind (see \eqref{eq39},\,\eqref{eq40}). Some recurrence relations for $\phi_{n,\lambda}^{(r)}(|z|^2)$ are obtained in Theorem 10. Another expression for $\phi_{n,\lambda}^{(r,r)}(|z|^2)$ is obtained in Theorem 11 by using the representation of the coherent state in terms of the number states. Finally, we define by introducing two new notations the unsigned degenerate Lah numbers and the signed degenerate Lah numbers, which are respectively degenerate versions of the Lah numbers and the signed Lah numbers. For the rest of this section, we recall the facts that are needed throughout this paper. \par

\medskip

For $n\geq0$, the Stirling numbers of the second kind are defined by
\begin{equation}\label{eq0}
\begin{split}
x^n=\sum_{k=0}^n S_2(n,k)(x)_k, \quad (n\geq0), \quad {(\text {see [1,4,5,10,14,16]})}.
\end{split}
\end{equation}

For any $\lambda \in \mathbb{R}$, the degenerate exponentials are given by
\begin{equation}\label{eq01}
\begin{split}
e_\lambda^x(t)=\sum_{k=0}^\infty (x)_{k,\lambda}\frac{t^k}{k!},\quad e_{\lambda}(t)=e_{\lambda}^{1}(t), \quad {(\text {see \cite{10,12}})},
\end{split}
\end{equation}
where the generalized falling factorials are given by
\begin{equation}\label{eq01-1}
\begin{split}
(x)_{0,\lambda}=1, (x)_{n,\lambda}=x(x-\lambda)\cdots(x-(n-1)\lambda), \quad (n\geq1).
\end{split}
\end{equation}

Recently, the degenerate Stirling numbers of the second kind are defined by
\begin{equation}\label{eq02}
\begin{split}
(x)_{n,\lambda}=\sum_{k=0}^n S_{2,\lambda}(n,k)(x)_k, \quad (n\geq0), \quad {(\text {see [7,8,10,11,12]})},
\end{split}
\end{equation}
where $(x)_0=1, (x)_n=x(x-1)\cdots(x-n+1), \quad (n\geq1)$. \\
Note that $\lim_{\lambda\rightarrow0}S_{2,\lambda}(n,k)=S_2(n,k), \quad (n, \ k \geq0)$.
From \eqref{eq02}, we note that
\begin{equation}\label{eq03}
\begin{split}
\frac{1}{k!}(e_\lambda(t)-1)^k=\sum_{n=k}^\infty S_{2,\lambda}(n,k)\frac{t^n}{n!}, \quad (k\geq0), \quad {(\text {see \cite{10}})}.
\end{split}
\end{equation}

It is well known that the ordinary Bell polynomials are defined by 
\begin{equation}\label{eq04}
\begin{split}
e^{x(e^t-1)}=\sum_{n=0}^\infty \phi_n(x)\frac{t^n}{n!}, \quad {(\text {see \cite{9, 13, 14}})}.
\end{split}
\end{equation}
When $x=1, \phi_n=\phi_n(1), \quad (n\geq0)$, are called the Bell numbers.\\
From \eqref{eq04}, we note that $\phi_n(x)=\sum_{k=0}^n S_2(n,k)x^k, \quad (n\geq0)$.

Recently, the degenerate Bell polynomials are given by 
\begin{equation}\label{eq05}
\begin{split}
e^{x(e_\lambda(t)-1)}=\sum_{n=0}^\infty \phi_{n,\lambda}(x)\frac{t^n}{n!}, \quad {(\text {see \cite{8,9,11,12}})}.
\end{split}
\end{equation}
By \eqref{eq05}, we get 
\begin{equation}\label{eq06}
\begin{split}
\phi_{n,\lambda}(x)=\sum_{k=0}^n S_{2,\lambda}(n,k)x^k, \quad (n\geq0).
\end{split}
\end{equation}
When $x=1, \phi_{n,\lambda}=\phi_{n,\lambda}(1), \quad (n\geq0)$, are called the degenerate Bell numbers. \par
Recall that $a$ and $a^\dag$ are the boson annihilation and creation operators such that
\begin{equation}\label{eq07}
\begin{split}
[a, a^\dag]=aa^\dag-a^\dag a=1, \quad {(\text {see [2,3,6-8,11-13,15]})}.
\end{split}
\end{equation}
The number states $|m\rangle, m=0, 1, \cdots,$ are given by
\begin{equation}\label{eq08}
\begin{split}
a|m\rangle=\sqrt{m}|m-1\rangle, \ a^\dag |m\rangle=\sqrt{m+1}|m+1\rangle.
\end{split}
\end{equation}
The coherent state $|z\rangle$, where $z$ is complex number, satisfies $a|z\rangle=z|z\rangle, \langle z|z \rangle=1$.
To show a connection to coherent states, we recall that the harmonic oscillator has Hamiltonian $H=a^\dag a$ (neglecting the zero point energy) and the usual eigenstates $|n\rangle (n\in \mathbb{N})$ satisfying
\begin{equation}\label{eq09}
\begin{split}
H|n\rangle=n|n\rangle \ \ {\rm{and}} \ \ \langle m|n\rangle=\delta_{m,n}, \quad \ \ {(\text {see \cite{8,11,13}})},
\end{split}
\end{equation}
where $\delta_{m,n} $ is Kronecker's symbol.


The normal ordering of a degenerate integral power of the number operator $a^\dag a$ in terms of the boson operators $a$ and $a^\dag$ can be written in the form
\begin{equation}\label{eq10}
\begin{split}
(a^\dag a)_{n,\lambda}=\sum_{k=0}^n S_{2,\lambda}(n,k)(a^\dag)^ka^k, \quad (n\geq0), \quad \ \ {(\text {see \cite{7,8,11}})}.
\end{split}
\end{equation}
We note that the standard bosonic commutation relation $[a,a^\dag]=aa^\dag-a^\dag a=1$ can be considered formally, in a suitable space of functions $f$, by letting $a=\frac{d}{dx}$ and $a^\dag=x$ (the operator of multiplication by $x$).

By \eqref{eq10}, we get
\begin{equation*}
\begin{split}
\bigg(x\frac{d}{dx}\bigg)_{n,\lambda}f(x)=\sum_{k=0}^n S_{2,\lambda}(n,k)x^k(\frac{d}{dx})^kf(x), \quad {(\text {see \cite{7,8,11}})}.
\end{split}
\end{equation*}
From the definition of coherent states, we note that $a|z\rangle=z|z\rangle$, equivalently $\langle z|a^\dag=\langle z|\overline{z}$, where $z \in \mathbb{C}$ and $\overline{z}$ is the complex conjugate of $z$.

In \cite{2,3}, Blasiak-Person-Solomon considered the generalized Stirling numbers of the second kind $S_{r,s}(n,k)$ which are given by
\begin{equation}\label{eq11}
\begin{split}
((a^{\dag})^r a^s)^n=(a^{\dag})^{n(r-s)}\sum_{k=s}^{ns}S_{r,s}(n,k)(a^{\dag})^k a^k,
\end{split}
\end{equation}
where $r,s$ are positive integers with $r\geq s$ and $n$ is any positive integer. They also considered the polynomials, which we call the $(r,s)$-Bell polynomials, given by
\begin{equation}\label{eq12}
\begin{split}
\phi_{n}^{(r,s)}(x)=\sum_{k=s}^{ns} S_{r,s}(n,k)x^{k}.
\end{split}
\end{equation}
For $x=1$, $\phi_{n}^{(r,s)}=\phi_{n}^{(r,s)}(1)$ are called the $(r,s)$-Bell numbers.

\medskip

\section{Generalized degenerate Stirling numbers arising from degenerate boson normal ordering}

\medskip

In this section, unless otherwise stated, $r,s$ are positive integers with $r\geq s$ and  $n$ is any positive integer.
In light of \eqref{eq11}, we introduce the {\it{generalized degenerate $(r,s)$-Stirling numbers of the second kind}} arising from the degenerate boson normal ordering
\begin{equation}\label{eq13}
\begin{split}
\prod_{k=0}^{n-1}\bigg[x^{r-s}\bigg(x^s\bigg(\frac{d}{dx}\bigg)^s-k\lambda\bigg)\bigg]
&=\prod_{k=0}^{n-1}\bigg(x^{r}\bigg(\frac{d}{dx}\bigg)^s -k\lambda x^{r-s}\bigg)\\
&=x^{n(r-s)}\sum_{k=0}^{ns}S_\lambda^{(r,s)}(n,k)x^k\bigg(\frac{d}{dx}\bigg)^k.
\end{split}
\end{equation}

From \eqref{eq11} and \eqref{eq13}, we see that $\lim_{\lambda \rightarrow 0}S_{\lambda}^{(r,s)}(n,k)=S_{r,s}(n,k)$.
 
\begin{remark}
(a) In \eqref{eq13} and below, the product of operators is understood to be written in the order of the factors corresponding to 1, 2, $\dots$, to $n-1$, from the left to the right.\\
(b) If $r=s$, then we see from \eqref{eq11} and \eqref{eq13} that $S_{\lambda}^{(r,r)}(n,k)=0$, for $ 0 \le k < r$.
\end{remark}

\noindent {\bf{Example:}} Let $n=2,\,r=4,\, s=2$ in \eqref{eq13}. Then, with $x=X,\, \frac{d}{dx}=D$, we have:
\begin{equation*}
\begin{split}
X^4D^2(X^4D^2-\lambda X^2)=(x^4D^2)^2-\lambda X^4D^2X^2,
\end{split}
\end{equation*}
where
\begin{equation*}
\begin{split}
D^2X^2&=D(XD+1)X=(DX)^2+DX=(XD+1)^2+XD+1\\
&=(XD)^2+3XD+2=X(XD+1)D+3XD+2=X^2D^2+4XD+2.
\end{split}
\end{equation*}
Thus, from \eqref{eq11}, we have
\begin{equation*}
\begin{split}
X^4D^2(X^4D^2-\lambda X^2)&=(x^4D^2)^2-\lambda X^6D^2-4\lambda X^5D-2 \lambda X^4\\
&=X^4\sum_{k=2}^{4}S_{4,2}(2,k)X^kD^k-\lambda X^6D^2-4\lambda X^5D-2 \lambda X^4\\
&=X^4\Big(S_{4,2}(2,4)X^4D^4+S_{4,2}(2,3)X^3D^3+(S_{4,2}(2,2)-\lambda)X^2D^2-4 \lambda XD-2 \lambda \Big)\\
&=X^4\Big(X^4D^4+8X^3D^3+(12-\lambda)X^2D^2-4 \lambda XD-2 \lambda \Big).
\end{split}
\end{equation*}
Thus $S_{\lambda}^{(4,2)}(2,4)=1,\, S_{\lambda}^{(4,2)}(2,3)=8,\, S_{\lambda}^{(4,2)}(2,2)=12-\lambda,\, S_{\lambda}^{(4,2)}(2,1)=-4\lambda,\, S_{\lambda}^{(4,2)}(2,0)=-2\lambda$.

\medskip

By applying the operators in \eqref{eq13} to $x^{p}$ and letting $x=1$, for any positive integer $p$, we obtain \eqref{eq14} with $x$ replaced by $p$ and hence it holds as polynomials.

\begin{equation}\label{eq14}
\begin{split}
\prod_{k=1}^n[(x+(k-1)(r-s))_s-(n-k)\lambda]=\sum_{k=0}^{ns} S_\lambda^{(r,s)}(n,k)(x)_k.
\end{split}
\end{equation}
From \eqref{eq02} and \eqref{eq14}, we note that $S_\lambda^{(1,1)}(n,k)=S_{2,\lambda}(n,k), \ \ (n, \ k\geq0)$. 

In view of \eqref{eq06}, we define the {\it{generalized degenerate $(r,s)$-Bell polynomials}} given by
\begin{equation}\label{eq15}
\begin{split}
\phi_{n,\lambda}^{(r,s)}(x)=\sum_{k=0}^{ns}S_\lambda^{(r,s)}(n,k)x^k.
\end{split}
\end{equation}
When $x=1, \phi_{n,\lambda}^{(r,s)}=\phi_{n,\lambda}^{(r,s)}(1)$ are called the {\it{generalized degenerate $(r,s)$-Bell numbers}}. Here we see that $\lim_{\lambda \rightarrow 0}\phi_{n,\lambda}^{(r,s)}(x)=\phi_{n}^{(r,s)}(x)$. \par
We observe that
\begin{equation}\label{eq16}
\begin{split}
e^{-x}\sum_{k=0}^\infty \frac{1}{k!}&\bigg(\prod_{j=1}^n[(k+(j-1)(r-s))_s-(n-j)\lambda]\bigg)x^k \\
&=e^{-x}\sum_{k=0}^\infty \frac{1}{k!}\sum_{l=0}^{ns}S_\lambda^{(r,s)}(n,l)(k)_lx^k \\
&=\sum_{l=0}^{ns}S_\lambda^{(r,s)}(n,l)e^{-x}\sum_{k=0}^\infty \frac{(k)_l}{k!}x^k \\
&=\sum_{l=0}^{ns}S_\lambda^{(r,s)}(n,l)e^{-x}x^l\bigg(\frac{d}{dx}\bigg)^l e^x \\
&=\sum_{l=0}^{ns}S_\lambda^{(r,s)}(n,l)x^le^{-x}e^x=\phi_{n,\lambda}^{(r,s)}(x).
\end{split}
\end{equation}
Therefore, by \eqref{eq16}, we obtain the following theorem.


\begin{theorem}
For $n\geq1$ and $r \ge s \ge 1$, we have
\begin{equation}\label{eq17}
\begin{split}
\phi_{n,\lambda}^{(r,s)}(x)=\frac{1}{e^x}\sum_{k=0}^\infty \frac{1}{k!}\bigg(\prod_{j=1}^n[(k+(j-1)(r-s))_s-(n-j)\lambda]\bigg)x^k.
\end{split}
\end{equation}
In particular, for $x=1$, we get
\begin{equation}\label{eq18}
\begin{split}
\phi_{n,\lambda}^{(r,s)}=\frac{1}{e}\sum_{k=0}^\infty \frac{1}{k!}\bigg(\prod_{j=1}^n[(k+(j-1)(r-s))_s-(n-j)\lambda]\bigg).
\end{split}
\end{equation}
\end{theorem}

\medskip

From \eqref{eq17}, we observe that
\begin{equation}\label{eq19}
\begin{split}
\phi_{n,\lambda}^{(r,r)}(x)=\frac{1}{e^x}\sum_{k=0}^{\infty}\frac{((k)_{r})_{n,\lambda}}{k!}x^{k}.
\end{split}
\end{equation}
In particular, from \eqref{eq06}, \eqref{eq15} and \eqref{eq19},  we note that
\begin{equation}\label{eq20}
\begin{split}
\phi_{n,\lambda}^{(1,1)}(x)=\frac{1}{e^x}\sum_{k=0}^\infty \frac{(k)_{n,\lambda}}{k!}x^{k}=\phi_{n,\lambda}(x).
\end{split}
\end{equation}

\medskip

From \eqref{eq18}, we observe that
\begin{equation}\label{eq21}
\begin{split}
\phi_n^{(r,s)}&=\lim_{\lambda\rightarrow0}\phi_{n,\lambda}^{(r,s)}=\frac{1}{e}\sum_{k=0}^\infty \frac{1}{k!}\prod_{j=0}^{n-1}(k+j(r-s))_s \\
&=\frac{1}{e}\sum_{k=0}^\infty \frac{1}{k!}\prod_{j=0}^{n-1}\big(k+j(r-s)\big)\big(k-1+j(r-s)\big)\cdots \big(k-s+1+j(r-s)\big) \\
&=\frac{1}{e}\sum_{k=0}^\infty \frac{1}{k!}\prod_{j=0}^{n-1}(r-s)^s\bigg(j+\frac{k}{r-s}\bigg)\bigg(j+\frac{k-1}{r-s}\bigg)\cdots\bigg(j+\frac{k-s+1}{r-s}\bigg) \\
&=\frac{1}{e}\sum_{k=0}^\infty
\frac{(r-s)^{sn}}{k!}\prod_{j=0}^{n-1} \prod_{l=1}^{s}\bigg(j+\frac{k-l+1}{r-s} \bigg)
=\frac{1}{e}\sum_{k=0}^\infty
\frac{(r-s)^{sn}}{k!}\prod_{l=1}^{s}\bigg(n-1+\frac{k-l+1}{r-s} \bigg)_{n}\\
&=\frac{1}{e}\sum_{k=0}^\infty
\frac{(r-s)^{sn}}{k!}\prod_{l=1}^{s}\frac{\Gamma(n+\frac{k-l+1}{r-s})}{\Gamma(\frac{k-l+1}{r-s})}=\frac{(r-s)^{sn}}{e}\sum_{k=0}^\infty
\frac{1}{k!}\prod_{l=1}^{s}\frac{\Gamma(n+\frac{k-l+1}{r-s})}{\Gamma(\frac{k-l+1}{r-s})}.
\end{split}
\end{equation}
Thus we obtain the following alternative expression for $\phi_n^{(r,s)}$.

\begin{theorem}
For $r > s \ge 1$ and $n \ge 1$, we have the following expression:
\begin{equation}\label{eq22}
\begin{split}
\phi_n^{(r,s)}=\frac{(r-s)^{sn}}{e}\sum_{k=0}^\infty
\frac{1}{k!}\prod_{l=1}^{s}\frac{\Gamma(n+\frac{k-l+1}{r-s})}{\Gamma(\frac{k-l+1}{r-s})}.
\end{split}
\end{equation}
\end{theorem}


From \eqref{eq14} and \eqref{eq17}, we have
\begin{equation}\label{eq23}
\begin{split}
\sum_{k=0}^{ns}&S_\lambda^{(r,s)}(n,k)x^k=\phi_{n,\lambda}^{(r,s)}(x) \\
&=e^{-x}\sum_{p=0}^\infty \frac{1}{p!}\bigg(\prod_{j=1}^n [(p+(j-1)(r-s))_s-(n-j)\lambda]\bigg)x^p \\
&=\sum_{m=0}^\infty (-1)^m\frac{x^m}{m!} \sum_{p=0}^\infty \bigg(\prod_{j=1}^n[(p+(j-1)(r-s))_s-(n-j)\lambda]\bigg)\frac{x^p}{p!} \\
&=\sum_{k=0}^\infty \sum_{p=0}^k \frac{(-1)^{k-p}k!}{(k-p)!p!}\bigg(\prod_{j=1}^n[(p+(j-1)(r-s))_s-(n-j)\lambda]\bigg)\frac{x^k}{k!} \\
&=\sum_{k=0}^\infty \frac{(-1)^k}{k!}\sum_{p=0}^k(-1)^p \binom{k}{p}\bigg(\prod_{j=1}^n[(p+(j-1)(r-s))_s-(n-j)\lambda]\bigg)x^k.
\end{split}
\end{equation}
Thus, by \eqref{eq23}, we obtain the following theorem.

\begin{theorem}
For $r \geq s\geq1$ and $n\geq1$, we have
\begin{equation}\label{eq24}
\begin{split}
\frac{(-1)^k}{k!}\sum_{p=0}^k(-1)^p\binom{k}{p}\bigg(\prod_{j=1}^n[(p+(j-1)(r-s))_s&-(n-j)\lambda]\bigg)\\
&=\left\{
\begin{split}
&S_\lambda^{(r,s)}(n,k), \ \ {\rm{if}} \ 0 \leq k\leq ns, \\
&0, \quad {\rm{if}} \ k > ns.
\end{split}\right.
\end{split}
\end{equation}
\end{theorem}

\medskip

From \eqref{eq13} and \eqref{eq24}, we note that
\begin{equation}\label{eq25}
\begin{split}
S_\lambda^{(r,s)}(n,k)&=\frac{(-1)^k}{k!}\sum_{p=0}^k (-1)^p\binom{k}{p}\prod_{j=1}^n\big[(p+(j-1)(r-s))_s-(n-j)\lambda\big] \\
&=\frac{(-1)^k}{k!}\sum_{p=0}^k (-1)^p\binom{k}{p}\prod_{j=0}^{n-1}\bigg[x^r\bigg(\frac{d}{dx}\bigg)^s-j\lambda x^{r-s}\bigg]x^p\Big|_{x=1} \\
&=\frac{(-1)^k}{k!}\prod_{j=0}^{n-1}\bigg[x^r\bigg(\frac{d}{dx}\bigg)^s-j\lambda x^{r-s}\bigg]\sum_{p=0}^k(-1)^p\binom{k}{p}x^p\Big|_{x=1} \\
&=\frac{(-1)^k}{k!}\prod_{j=0}^{n-1}\bigg[\bigg(x^r\bigg(\frac{d}{dx}\bigg)^s-j\lambda x^{r-s}\bigg)\bigg](1-x)^k \Big|_{x=1}.
\end{split}
\end{equation}
Therefore, by \eqref{eq25}, we obtain the following theorem.


\begin{theorem}
For $r \geq s\geq1$ and $n\geq1$, we have
\begin{equation}\label{eq26}
\begin{split}
S_\lambda^{(r,s)}(n,k)=\frac{(-1)^k}{k!}\prod_{j=0}^{n-1}\bigg[\bigg(x^r\bigg(\frac{d}{dx}\bigg)^s-j\lambda x^{r-s}\bigg)\bigg](1-x)^k \Big|_{x=1}.
\end{split}
\end{equation}
\end{theorem}

\medskip

From \eqref{eq26}, we note that
\begin{equation*}
\begin{split}
S_{2,\lambda}(n,k)=S_\lambda^{(1,1)}(n,k)&=\frac{(-1)^k}{k!}\prod_{j=0}^{n-1}\bigg[\bigg(x\frac{d}{dx}-j\lambda\bigg)\bigg](1-x)^k\Big|_{x=1} \\
&=\frac{(-1)^k}{k!}\bigg(x\frac{d}{dx}\bigg)_{n,\lambda}\sum_{p=0}^k\binom{k}{p}(-1)^px^p\Big|_{x=1} \\
&=\frac{(-1)^k}{k!}\sum_{p=0}^k\binom{k}{p}(-1)^p(p)_{n,\lambda},\quad (n \ge 1).
\end{split}
\end{equation*}


From \eqref{eq14}, we note that
\begin{equation}\label{eq27}
\begin{split}
\sum_{k=0}^{nr}S_\lambda^{(r,r)}(n,k)(x)_k=\prod_{k=1}^n[(x)_r-(n-k)\lambda]=((x)_r)_{n,\lambda}.
\end{split}
\end{equation}

Thus, by \eqref{eq27}, we get
\begin{equation}\label{eq28}
\begin{split}
((x)_r)_{n,\lambda}=\sum_{k=0}^{nr}S_\lambda^{(r,r)}(n,k)(x)_k.
\end{split}
\end{equation}

In particular, for $r=1$, we have
\begin{equation}\label{eq29}
\begin{split}
\sum_{k=0}^nS_\lambda^{(1,1)}(n,k)(x)_k=(x)_{n,\lambda}=\sum_{k=0}^nS_{2,\lambda}(n,k)(x)_k.
\end{split}
\end{equation}

Thus, by \eqref{eq29}, we get
\begin{equation}\label{eq30}
\begin{split}
S_\lambda^{(1,1)}(n,k)=S_{2,\lambda}(n,k), \  \ (n, \ k\geq0).
\end{split}
\end{equation}

From \eqref{eq26}, we note that
\begin{equation}\label{eq31}
\begin{split}
S_\lambda^{(r,r)}(n,k)&=\frac{(-1)^k}{k!}\prod_{j=0}^{n-1}\bigg(x^r\bigg(\frac{d}{dx}\bigg)^r-j\lambda\bigg)\sum_{p=0}^k(-1)^p\binom{k}{p}x^p\Big|_{x=1} \\
&=\frac{(-1)^k}{k!}\sum_{p=0}^k(-1)^p\binom{k}{p}((p)_r)_{n,\lambda}.
\end{split}
\end{equation}
Therefore, by \eqref{eq31}, we obtain the following theorem.


\begin{theorem}
For $r,\, n \geq1$, we have
\begin{equation}\label{eq32}
\begin{split}
S_\lambda^{(r,r)}(n,k)=\frac{(-1)^k}{k!}\sum_{p=0}^k(-1)^p\binom{k}{p}((p)_r)_{n,\lambda}.
\end{split}
\end{equation}
\end{theorem}

\medskip

From \eqref{eq32}, we note that
\begin{equation}\label{eq33}
\begin{split}
S_\lambda^{(r,r)}(1,r)=\frac{(-1)^r}{r!}\sum_{p=r}^r\binom{r}{p}(p)_r(-1)^p=\frac{(-1)^r}{r!}\binom{r}{r}(-1)^rr!=1.\\
\end{split}
\end{equation}

We recall that $\phi_{n,\lambda}$ and $S_{2,\lambda}(n,k)$ are related to special quantum states, called coherent states, defined as linear combinations of the eigenstates of the harmonic oscillator, $H=a^\dag a$, $H|n\rangle=n|n\rangle$, $\langle n|m\rangle=\delta_{m,n}$ and defined as $|z\rangle=e^{-\frac{|z|^2}{2}}\sum_{n=0}^\infty \frac{z^n}{\sqrt{n}}|n\rangle$,  with $\langle z|z\rangle=1$, for complex $z$,\\
 (see [3,7,11-14]).


From \eqref{eq13}, we have
\begin{equation}\label{eq34}
\begin{split}
\prod_{k=0}^{n-1}[(a^\dag)^{r-s}((a^\dag)^sa^s-k\lambda)] 
=(a^\dag)^{n(r-s)}\sum_{k=0}^{rs}S_\lambda^{(r,s)}(n,k)(a^\dag)^ka^k.
\end{split}
\end{equation}

\medskip

By \eqref{eq20} and \eqref{eq34}, we get
\begin{equation}\label{eq35}
\begin{split}
\langle z|e_\lambda^{a^\dag a}(t)|z\rangle &=\sum_{n=0}^\infty\frac{t^n}{n!} \langle z|(a^\dag a)_{n,\lambda}|z\rangle \\
&=\sum_{n=0}^\infty \frac{t^n}{n!}\langle z|\prod_{k=1}^n(a^\dag a-(n-k)\lambda)|z\rangle \\
&=\sum_{n=0}^\infty \frac{t^n}{n!}\sum_{k=0}^n S_{2,\lambda}^{(1,1)}(n,k)(\overline{z})^kz^k\langle z|z\rangle \\
&=\sum_{n=0}^\infty \frac{t^n}{n!}\phi_{n,\lambda}^{(1,1)}(|z|^2)=e^{|z|^2 (e_\lambda(t)-1)}.
\end{split}
\end{equation}


From \eqref{eq31} and \eqref{eq34}, we note that
\begin{equation}\label{eq36}
\begin{split}
\langle z|\prod_{k=0}^{n-1}[(a^\dag)^r a^r-k\lambda]|z\rangle&=\sum_{k=0}^{nr}S_\lambda^{(r,r)}(n,k)\langle z|(a^\dag)^ka^k|z\rangle 
=\sum_{k=0}^{nr}S_\lambda^{(r,r)}(n,k)(\overline{z})^kz^k\langle z|z\rangle\\
&=\sum_{k=0}^{nr}S_\lambda^{(r,r)}(n,k)(|z|^2)^k =\phi_{n,\lambda}^{(r,r)}(|z|^2)\\
&=\sum_{k=0}^{nr}\frac{(-1)^k}{k!}\sum_{p=0}^k(-1)^p\binom{k}{p}((p)_r)_{n,\lambda}(|z|^2)^{k}\\
&=\sum_{p=0}^{nr}\sum_{k=p}^{nr}\frac{(-1)^{k-p}}{k!}\binom{k}{p}((p)_r)_{n,\lambda}(|z|^2)^{k}.
\end{split}
\end{equation}

In particular, when  $|z|=1$, we have
\begin{equation}\label{eq37}
\begin{split}
\langle z|\prod_{k=0}^{n-1}[(a^\dag)^r a^r-k\lambda]|z\rangle&=\phi_{n,\lambda}^{(r,r)}=\sum_{k=0}^{nr}S_\lambda^{(r,r)}(n,k)\\
&=\sum_{p=0}^{nr}\sum_{k=p}^{nr}\frac{(-1)^{k-p}}{k!}\binom{k}{p}((p)_r)_{n,\lambda}.
\end{split}
\end{equation}
Therefore, by \eqref{eq36} and \eqref{eq37}, we obtain the following theorem.

\begin{theorem}
For $r,\,n\geq1$, we have
\begin{equation*}
\begin{split}
\phi_{n,\lambda}^{(r,r)}(|z|^2)=\langle z|\prod_{k=0}^{n-1}[(a^\dag)^ra^r-k\lambda]|z\rangle=\sum_{p=0}^{nr}\sum_{k=p}^{nr}\frac{(-1)^{k-p}}{k!}\binom{k}{p}((p)_r)_{n,\lambda}(|z|^2)^{k}.
\end{split}
\end{equation*}

In particular, when  $|z|=1$, we get
\begin{equation*}
\begin{split}
\phi_{n,\lambda}^{(r,r)}=\langle z|\prod_{k=0}^{n-1}[(a^\dag)^ra^r-k\lambda]|z\rangle=\sum_{p=0}^{nr}\sum_{k=p}^{nr}\frac{(-1)^{k-p}}{k!}\binom{k}{p}((p)_r)_{n,\lambda}.
\end{split}
\end{equation*}
\end{theorem}

\medskip


Now, we observe  \eqref{eq36} that
\begin{equation}\label{eq38}
\begin{split}
\langle z|e_\lambda^{(a^\dag)^ra^r}(t)|z\rangle&=\sum_{n=0}^\infty \frac{t^n}{n!}\langle z|((a^\dag)^ra^r)_{n,\lambda}|z\rangle \\
&=\sum_{n=0}^\infty \frac{t^n}{n!}\langle z|\prod_{k=0}^{n-1}[(a^\dag)^ra^r-k\lambda]|z\rangle \\
&=\sum_{n=0}^\infty \frac{t^n}{n!}\phi_{n,\lambda}^{(r,r)}(|z|^2).
\end{split}
\end{equation}
Therefore, by \eqref{eq38}, we obtain the following theorem.

\begin{theorem}
Let $r$ be a positive integer.  Then the generating function of $\phi_{n,\lambda}^{(r,r)}(|z|^2)$ is given by
\begin{equation*}
\begin{split}
\langle z|e_\lambda^{(a^\dag)^ra^r}(t)|z\rangle=\sum_{n=0}^\infty \phi_{n,\lambda}^{(r,r)}(|z|^2)\frac{t^n}{n!}.
\end{split}
\end{equation*}
\end{theorem}

\medskip

We recall that the degenerate $r$-Stirling numbers of the second kind are define by
\begin{equation}\label{eq39}
\begin{split}
(x+r)_{n,\lambda}=\sum_{k=0}^n {n+r \brace k+r}_{r,\lambda}(x)_k, \ \ (n\geq0), \ \ {(\text {see \cite{8,12}})}.
\end{split}
\end{equation}

The degenerate $r$-Bell polynomials are given by
\begin{equation}\label{eq40}
\begin{split}
\phi_{n,\lambda}^{(r)}(x)=\sum_{k=0}^n {n+r \brace k+r}_{r,\lambda} x^k, \ \ (n\geq0), \ \ {(\text {see \cite{8,12}})}.
\end{split}
\end{equation}
For $x=1$, $\phi_{n,\lambda}^{(r)}=\phi_{n,\lambda}^{(r)}(1)$ are called the degenerated $r$-Bell numbers.

Now, we recall from \cite{8} that
\begin{equation}\label{eq41}
\begin{split}
(a^{\dag}a+r)_{n,\lambda}=\sum_{k=0}^n{n+r \brace k+r}_{r,\lambda} (a^\dag)^ka^k, \ \ (n\geq0).
\end{split}
\end{equation}

From \eqref{eq41}, we have
\begin{equation}\label{eq42}
\begin{split}
\langle z|e_\lambda^{a^{\dag}a+r}(t)|z\rangle&=\sum_{n=0}^\infty \frac{t^n}{n!}\langle z|(a^{\dag}a+r)_{n,\lambda}|z\rangle \\
&=\sum_{n=0}^\infty \frac{t^n}{n!}\sum_{k=0}^n {n+r \brace k+r}_{r,\lambda}\langle z|(a^\dag)^ka^k|z\rangle \\
&=\sum_{n=0}^\infty \frac{t^n}{n!}\sum_{k=0}^n {n+r \brace k+r}_{r,\lambda}(\overline{z})^kz^k\langle z|z\rangle \\
&=\sum_{n=0}^\infty \frac{t^n}{n!}\sum_{k=0}^n {n+r \brace k+r}_{r,\lambda}(|z|^2)^k=\sum_{n=0}^\infty \frac{t^n}{n!}\phi_{n,\lambda}^{(r)}(|z|^2) \\
&=e_\lambda^r(t)e^{|z|^2(e_\lambda(t)-1)}.
\end{split}
\end{equation}
Therefore, by \eqref{eq42}, we obtain the following theorem.

\begin{theorem}
Let $r$ be nonnegative integer. Then we have
\begin{equation*}
\begin{split}
\sum_{n=0}^\infty \phi_{n,\lambda}^{(r)}(|z|^2)\frac{t^n}{n!}=\langle z|e_\lambda^{a^{\dag}a+r}(t)|z\rangle=e_\lambda^r(t)e^{|z|^2(e_\lambda(t)-1)}.
\end{split}
\end{equation*}
\end{theorem}

\medskip

\medskip

Let $g(t)=\sum_{n=0}^\infty \phi_{n,\lambda}^{(r)}(|z|^2)\frac{t^n}{n!}=\langle z| e_{\lambda}^{a^{\dag}a+r}(t)|z \rangle=e_\lambda^r(t)e^{|z|^2(e_\lambda(t)-1)}$. Then from the series expression we have 
\begin{equation}\label{eq43}
\begin{split}
\frac{dg(t)}{dt}=\sum_{n=0}^{\infty}\phi_{n+1,\lambda}^{(r)}(|z|^2)\frac{t^n}{n!}.
\end{split}
\end{equation}
From the expression of $g(t)$ in the bracket notation, we get
\begin{equation}\label{eq44}
\begin{split}
\frac{dg(t)}{dt}&=\langle z| (a^{\dag}a+r)e_{\lambda}^{a^{\dag}a+r-\lambda}(t)|z \rangle \\
&=e_{\lambda}^{-\lambda}(t)\Big(\langle z| a^{\dag} e_{\lambda}^{aa^{\dag}+r}(t) a|z \rangle + r \langle z| e_{\lambda}^{a^{\dag}a+r}(t)|z \rangle \Big)\\
&=e_{\lambda}^{-\lambda}(t)\Big(|z|^2\langle z| e_{\lambda}^{a^{\dag}a+r+1}(t) |z \rangle + r \langle z| e_{\lambda}^{a^{\dag}a+r}(t)|z \rangle \Big)\\
&=e_{\lambda}^{-\lambda}(t)\Big(|z|^2\sum_{k=0}^\infty \phi_{k,\lambda}^{(r+1)}(|z|^2)\frac{t^k}{k!} + r \sum_{k=0}^\infty \phi_{k,\lambda}^{(r)}(|z|^2)\frac{t^k}{k!} \Big)\\
&=\sum_{l=0}^{\infty}(-\lambda)^l l! \frac{t^l}{l!}\sum_{k=0}^{\infty}\Big(|z|^2\phi_{k,\lambda}^{(r+1)}(|z|^2)+r \phi_{k,\lambda}^{(r)}(|z|^2) \Big)\frac{t^k}{k!}\\
&=\sum_{n=0}^{\infty}\bigg(\sum_{k=0}^{n}\binom{n}{k}(-\lambda)^{n-k}(n-k)!\Big(|z|^2\phi_{k,\lambda}^{(r+1)}(|z|^2)+r\phi_{k,\lambda}^{(r)}(|z|^2)\Big)\bigg)\frac{t^n}{n!}.
\end{split}
\end{equation}
From the last expression of $g(t)$, we obtain
\begin{equation}\label{eq45}
\begin{split}
\frac{dg(t)}{dt}&=\Big(re_{\lambda}^{-\lambda}(t)+|z|^2e_{\lambda}^{1-\lambda}(t)\Big)e_{\lambda}^{r}(t)e^{|z|^2(e_{\lambda}(t)-1)}\\
&=\sum_{k=0}^{\infty}\big(r(-\lambda)_{k,\lambda}+|z|^2(1-\lambda)_{k,\lambda}\big)\frac{t^k}{k!}\sum_{m=0}^{\infty}\phi_{m,\lambda}^{(r)}(|z|^2)\frac{t^m}{m!}\\
&=\sum_{n=0}^{\infty}\sum_{k=0}^{n}\binom{n}{k}\big(r(-\lambda)_{k,\lambda}+|z|^2(1-\lambda)_{k,\lambda}\big)\phi_{n-k,\lambda}^{(r)}(|z|^2)\frac{t^n}{n!}.
\end{split}
\end{equation}

Now, from \eqref{eq43}, \eqref{eq44} and \eqref{eq45}, we obtain the next result.

\begin{theorem}
Let $r,\,n$ be nonnegative integers. Then we have the following recurrence relations:
\begin{equation*}
\begin{split}
\phi_{n+1,\lambda}^{(r)}(|z|^2)&=\sum_{k=0}^{n}\binom{n}{k}(-\lambda)^{n-k}(n-k)!\Big(|z|^2\phi_{k,\lambda}^{(r+1)}(|z|^2)+r\phi_{k,\lambda}^{(r)}(|z|^2)\Big)\\
&=\sum_{k=0}^{n}\binom{n}{k}\big(r(-\lambda)_{k,\lambda}+|z|^2(1-\lambda)_{k,\lambda}\big)\phi_{n-k,\lambda}^{(r)}(|z|^2).
\end{split}
\end{equation*}
In particular, for $|z|=1$, we get
\begin{equation*}
\begin{split}
\phi_{n+1,\lambda}^{(r)}&=\sum_{k=0}^{n}\binom{n}{k}(-\lambda)^{n-k}(n-k)!\Big(\phi_{k,\lambda}^{(r+1)}+r\phi_{k,\lambda}^{(r)}\Big)\\
&=\sum_{k=0}^{n}\binom{n}{k}\big(r(-\lambda)_{k,\lambda}+(1-\lambda)_{k,\lambda}\big)\phi_{n-k,\lambda}^{(r)}.
\end{split}
\end{equation*}
\end{theorem}

\medskip

Evaluating the left hand side of \eqref{eq36} by using the representation of the coherent state in terms of the number states, we have
\begin{equation}\label{eq46}
\begin{split}
&\langle z|((a^\dag)^ra^r)_{k,\lambda}|z\rangle=\langle z|\Pi_{l=0}^{k-1}[(a^\dag)^ra^r-l\lambda]|z\rangle \\
&=e^{-\frac{|z|^2}{2}}e^{-\frac{|z|^2}{2}}\sum_{m,n=0}^\infty \frac{(\overline{z})^m(z)^n}{\sqrt{m!}\sqrt{n!}}((n)_r)_{k,\lambda}\langle m|n\rangle \\
&=e^{-|z|^2}\sum_{n=0}^\infty \frac{(|z|^2)^n}{n!}((n)_r)_{k,\lambda}=e^{-|z|^2}\sum_{n=1}^\infty \frac{(|z|^2)^n}{n!}((n)_r)_{k,\lambda},
\end{split}
\end{equation}
where $r$ and $k$ are positive integer.\par
Thus, from \eqref{eq36} and \eqref{eq46}, we get the next theorem.
\begin{theorem}
Let $r, \, k$ be positive integers. Then we have
\begin{equation*}\label{eq47}
\begin{split}
\phi_{k,\lambda}^{(r,r)}(|z|^2)=e^{-|z|^2}\sum_{n=1}^\infty \frac{(|z|^2)^n}{n!}((n)_r)_{k,\lambda}.
\end{split}
\end{equation*}
 In particular, when $|z|=1$, we get
\begin{equation*}\label{eq48}
\begin{split}
\phi_{k,\lambda}^{(r,r)}=\frac{1}{e}\sum_{n=1}^\infty \frac{1}{n!}((n)_r)_{k,\lambda}.
\end{split}
\end{equation*}
\end{theorem}

Lastly, we introduce two notations given by
\begin{equation}\label{eq49}
\begin{split}
\langle \langle x \rangle \rangle_{0,\lambda}=1, \langle \langle x \rangle \rangle_{n,\lambda}=\prod_{k=1}^n (x+(k-1)-(n-k)\lambda), \quad (n \geq 1).
\end{split}
\end{equation}

and

\begin{equation}\label{eq50}
\begin{split}
((x))_{0,\lambda}=1, ((x))_{n,\lambda}=\prod_{k=1}^n (x-(k-1)+(n-k)\lambda), \quad (n \geq 1).
\end{split}
\end{equation}

We may consider the unsigned degenerate Lah numbers defined by
\begin{equation}\label{eq51}
\begin{split}
\langle \langle x \rangle \rangle_{n,\lambda}=\sum_{k=0}^n L_\lambda(n,k)(x)_k, \quad (k\geq0), \quad (n\geq0).
\end{split}
\end{equation}

In addition, the signed degenerate Lah numbers are given by
\begin{equation}\label{eq52}
\begin{split}
((x))_{n,\lambda}=\sum_{k=0}^n L_\lambda^1(n,k)\langle x \rangle_k, \quad (n\geq0).
\end{split}
\end{equation}
Note that $\lim_{\lambda \rightarrow 0} L_\lambda(n,k)=L(n,k), \quad (n, \ k\geq0)$ are the ordinary Lah numbers given by $\langle x \rangle_n=\sum_{k=0}^nL(n,k)(x)_k$, where $\langle x \rangle_0=1, \langle x \rangle_n=x(x+1)\cdots(x+n-1), \quad (n\geq1)$, \ \ (see \cite{3}). Moreover, $\lim_{\lambda \rightarrow 0}L_{\lambda}^{1}(n,k)=L^{1}(n,k)$, where $L^{1}(n,k)$ are the signed Lah numbers given by $(x)_{n}=\sum_{k=0}^{n}L^{1}(n,k)\langle x \rangle_{k}$. We also observe from \eqref{eq14} and \eqref{eq51} that
\begin{equation*}
\begin{split}
S_\lambda^{(2,1)}(n,k)=L_\lambda(n,k), \quad (n \ge 1).
\end{split}
\end{equation*}

\section{Conclusion}
In recent years, studying degenerate versions of some special numbers and polynomials has drawn the attention of many mathematicians and yielded many interesting results. These degenerate versions include the degenerate Stirling numbers of the first and second kinds, degenerate Bernoulli numbers of the second kind and degenerate Bell numbers and polynomials. \par
In this paper, we introduced the generalized degenerate $(r,s)$-Stirling numbers of the second kind and the generalized degenerate $(r,s)$-Bell polynomials by considering the boson normal ordering of $ \prod_{k=0}^{n-1}\big((a^{\dag})^{r}a^s -k\lambda (a^{\dag})^{r-s}\big)$. We studied some properties, explicit expressions and generating functions of those numbers and polynomials. They are degenerate versions of the corresponding ones in the earlier works by Blasiak-Person-Solomon (see \cite{4}). \par
These new numbers are expected to play an important role in the study of various degenerate versions of many special polynomials and numbers, just as the degenerate Stirling numbers of the second kind have played an important role in that study.\par
We would like to continue to study various degenerate versions of many special polynomials and numbers and their applications to physics, science and engineering as well as to mathematics. 

\bigskip

\noindent{\bf{Acknowledgments}} \\
The authors thank Jangjeon Institute for Mathematical Science for the support of this research.

\vspace{0.1in}

\noindent{\bf {Availability of data and material}} \\
Not applicable.

\vspace{0.1in}

\noindent{\bf{Funding}} \\
The third author is  supported by the Basic Science Research Program, the National
               Research Foundation of Korea,
               (NRF-2021R1F1A1050151).
\vspace{0.1in}

\noindent{\bf{Ethics approval and consent to participate}} \\
The authors declare that there is no ethical problem in the production of this paper.

%
%
%
%

%

\

\end{document}